\documentclass{amsart}
\usepackage{hyperref}
\usepackage{algpseudocode}
\usepackage{amssymb}
\usepackage{amsmath}

\def\bne{\begin{equation}}
\def\ene{\end{equation}}

\def\Re{\mathbb{R}}

\newtheorem{algorithm}{Algorithm}

\hypersetup{
    pdftitle={A Fast Compensated Algorithm for Computing Givens Rotations},
    pdfauthor={Carlos F. Borges},
    pdfkeywords={Givens Rotation, compensated, fused multiply-add, IEEE 754}
    }

\def\Re{\mathbb{R}}

\title{A Fast Compensated Algorithm for Computing Givens Rotations}

\author{Carlos F. Borges}
\address{Department of Applied Mathematics\\Naval Postgraduate School\\Monterey CA 93943\\borges@nps.edu}

\begin{document}

\begin{abstract}
We develop a very simple compensated scheme for computing very accurate Givens rotations. The approach is significantly more straightforward than the one in \cite{borges2021fast}, and the derivation leads to a very satisfying algorithm whereby a naively computed Givens rotation can be used to construct a correction to itself. It is also seen that this scheme continues to provide high accuracy even when built on a hypoteneuse calculation that is of lesser accuracy.
\end{abstract}

\maketitle

\section{Mathematical Preliminaries}

Givens rotations \cite{Wilkinson65,GVL,Demmel1997} are a fundamental tool in numerical linear algebra. Given two real numbers $f,g \in \Re$ we define a real Givens rotation $R(c,s) \in \Re^{2\times2}$ to be a real unitary matrix such that
$$
R(c,s) \cdot \begin{bmatrix}f\\g\end{bmatrix} \equiv \begin{bmatrix}
    c & s \\
    -s & c
\end{bmatrix} \cdot \begin{bmatrix}f\\g\end{bmatrix} = \begin{bmatrix}r\\0\end{bmatrix}.
$$
Since $R(c,s)$ is real and unitary it is clear that
\begin{equation}
   c^2 + s^2 = 1.
   \label{normality}
\end{equation}
The conventional approach is to take
\begin{eqnarray*}
   r & = & \sqrt{f^2+g^2} \\
   c & = & \frac{f}{r} \\
   s & = & \frac{g}{r}
\end{eqnarray*}
provided that $r > 0$. There are several possible conventions for constructing Givens rotations and their varying strengths and weaknesses are discussed in \cite{Anderson,BindelDemmelKahan} and others. We will not debate those here but will instead rely on the convention that is embodied in the DLARTG\footnote{This is the double precision version but all that follows will also apply to the single precision SLARTG code provided that the machine constants are correct for the floating point format that is chosen.} code from LAPACK version 3.12.0 which appears in appendix A.

First observe that the code deals with the possibility that $r = 0$ with some lead in branching and includes protections against avoidable overflow/underflow. We make two observations. First is that the second branch, from lines 7 to 10, is superfluous. Removing it will not change the results in any way. Second, the final branch, from lines 17 to 25, is there to prevent avoidable overflow/underflow in the computation of $r$. That means this can also be removed if the programming environment provides a proper {\tt hypot(f,g)} utility, as recommended in the IEEE754 standard, by replacing the right hand side of the assignment on line 13 with a call to such a utility.\footnote{This aspect of the change will be of particular interest to us as we will see that there can be advantages to different approaches to the {\tt hypot()} calculation.}

In sum, invoking these simplifications we proceed using the following algorithm:

\begin{algorithm} Simplified DLARTG

   \hrule
   \begin{algorithmic}
      \If{$g == 0$}
         \State{$c = 1$}
         \State{$s = 0$}
         \State{$r = f$}
      \Else
      \State{$d = hypot(f,g)$}
      \State{$c = abs(f) / d$}
      \State{$r = copysign( d, f )$}
      \State{$s = g / r$}
      \EndIf
   \end{algorithmic}
   \hrule
   \label{DLARTG}
\end{algorithm}

\section{The influence of the {\tt hypot()} calculation}

The {\tt hypot()} calculation is a critical part of constructing Givens rotations and we have seen that it strongly influences the structure or the LAPACK DLARTG code. The critical issue is that the calculation can lead to avoidable overflow/underflow errors if not done carefully. We have noted that using a {\tt hypot()} utility simplifies the code but wonder how different approaches to the {\tt hypot()} calculation affect the construction of a Givens rotation. We note that there are three main approaches. First is to use a correctly rounded {\tt hypot()} utility function such as the one developed in \cite{borges}. Second, one can use the naive calculation $r = \sqrt{f^2*g^2}$ paired with scaling by powers of the radix to prevent avoidable overflow/underflow as in the LAPACK code. This approach is very accurate but often does not yield correctly rounded results. And third, one can use the following common device\footnote{This is a common device but we note that this particular implementation assumes that at least one of the inputs is non-zero.} for computing $r$ without avoidable overflow/underflow:

 \vspace{.2in}
\begin{algorithm} WeakHypot(f,g)
\hrule
\begin{algorithmic}
   \State $af = abs(f)$
   \State $ag = abs(g)$
   \If{$ag > af$}
       \State $af, ag = ag, af$
   \EndIf
   \State $r = ag / af$
   \State $r = af * sqrt(1 + r * r)$
\end{algorithmic}
\hrule
\label{WeakHypot}
\end{algorithm}
\vspace{.2in}

This approach is very fast but less accurate than the naive approach (see \cite{borges}).

To understand how the accuracy of the {\tt hypot()} calculation influences the accuracy of the Givens rotation we perform a simple test. We begin by computing values for $c$ and $s$ using algorithm \ref{DLARTG} with variables represented in the BigFloat format in Julia and then round the results to get our baseline values. We note that this computation cannot be guaranteed to generate correctly rounded values for $c$ and $s$ (although it can be expected to do so in nearly every case). We then compare the outputs of algorithm \ref{DLARTG} operating in Float64 when using three different algorithms for the {\tt hypot()} calculation against the baseline values to see how well they match. We will use $10^{9}$ normally distributed random inputs, that is both $x,y \sim \mathcal{N}(0,1)$. In the table below we summarize the error rates.

\begin{table}[h]
   \begin{center}
       \caption{Error rates for computing the Givens rotation with the simplified DLARTG algorithm using different variations of the {\tt hypot()} calculation}
       \label{tab:tableGivens2}
       \begin{tabular}{l|l|l|l|l|l|l}
          & \multicolumn{2}{c}{Correct} & \multicolumn{2}{c}{Naive} & \multicolumn{2}{c}{Weak}\\
          & Cosine & Sine & Cosine & Sine & Cosine & Sine \\
          \hline
          Zero ulp errors & 71.069 & 71.062 & 66.563 & 66.567 & 51.700 & 54.700 \\
          One ulp errors  & 28.940 & 28.938 & 33.207 & 33.204 & 45.200 & 42.800 \\
          Two ulp errors  & 0      & 0      &  0.230 &  0.230 &  3.100 &  2.500
       \end{tabular}
    \end{center}
 \end{table}

In no case do the results match up perfectly with the higher precision calculation. However, we can now clearly see the trade-off between the various approaches - the more accurate (and costly) the {\tt hypot()} calculation, the more accurate the computed Givens rotation.

\section{A Compensated Approach to Givens Rotations}

Although the true quantities $c$ and $s$ satisfy four mathematical conditions -
\begin{eqnarray}
   cs - sc & = & 0 \label{orth1}\\
   c^2 + s^2 & = & 1 \label{norm1}\\
   cg - sf & = & 0 \label{orth2}\\
   cf + sg & = & r \label{norm2}
\end{eqnarray}
the computed quantities $\bar{c}$ and $\bar{s}$ may not. In fact, with the exception of condition \ref{orth1}, which is satisfied by the computed quantities in any IEEE754 compliant FPS since multiplication is commutative, the remaining three conditions are frequently not satisfied by the computed values, no matter how accurately they are computed.

We examine the problem by writing the true values of $c$ and $s$ as simple additive perturbations to the computed values by setting
\begin{eqnarray*}
   c & = & \bar{c} + \delta_c \\
   s & = & \bar{s} + \delta_s
\end{eqnarray*}
Conditions \ref{norm1}, \ref{orth2}, \ref{norm2} become
\begin{eqnarray*}
   (\bar{c} + \delta_c)^2 + (\bar{s} + \delta_s)^2 & = & 1 \\
   (\bar{c} + \delta_c)g - (\bar{s}+\delta_s)f & = & 0 \\
   (\bar{c} + \delta_c)f + (\bar{s} + \delta_s)g & = & r
\end{eqnarray*}
which can be rearranged to the following set of three equations for the two unknown perturbations $\delta_c$ and $\delta_s$
\begin{eqnarray}
   2 \bar{c} \delta_c + 2 \bar{s} \delta_s & = & 1 - \bar{c}^2 - \bar{s}^2 - \delta_c^2 - \delta_s^2 \label{norm1p}\\
   -g \delta_c + f \delta_s & = & \bar{c}g - \bar{s}f \label{orth2p}\\
   f \delta_c + g \delta_s & = & r - \bar{c}f - \bar{s}g \label{norm2p}.
\end{eqnarray}

Experience indicates that the most useful approach to solving these equations, at least approximately, is to ignore the tiny quantities $\delta_{c}^2$ and $\delta_{s}^2$ in equation \ref{norm1p} and to drop equation \ref{norm2p} from consideration.\footnote{We note that equation \ref{norm2p} is far less important than equation \ref{norm1p} as the latter seeks to guarantee that the compensated Givens rotation preserves the norm of all vectors, whereas the former only seeks to guarantee that it preserves the norm of one specific vector.} Furthermore, dividing equation \ref{orth2p} by $r$ and replacing $f/r$ and $g/r$ by the computed values $\bar{c}$ and $\bar{s}$ gives
\begin{equation}
   \begin{bmatrix}
      \bar{c} & \bar{s} \\ -\bar{s} & \bar{c}
   \end{bmatrix}
   \begin{bmatrix}
      \delta_c \\ \delta_s
   \end{bmatrix} =
   \begin{bmatrix}
      (1 - \bar{c}^2 - \bar{s}^2)/2 \\ (\bar{c}g - \bar{s}f)/r
   \end{bmatrix}
   \label{BasicEq}
\end{equation}
And this is solved, at least approximately, by simply applying the transpose of the originally computed Givens rotation. That is
\begin{equation}
   \begin{bmatrix}
      \delta_c \\ \delta_s
   \end{bmatrix} = \begin{bmatrix}
      \bar{c} & -\bar{s} \\ \bar{s} & \bar{c}
   \end{bmatrix}
   \begin{bmatrix}
      (1 - \bar{c}^2 - \bar{s}^2)/2 \\ (\bar{c}g - \bar{s}f)/r
   \end{bmatrix}
\end{equation}

It is critical to compute the quantities on the right hand side of \ref{BasicEq} as accurately as possible since they both involve {\em significant} cancellation. To do this we can use three applications of a double length product such as the {\tt 2MultFMA()} algorithm (see \cite{FPHB}) and one additional {\tt fma()} call. In particular, the {\tt 2MultFMA()} algorithm represents the product to two floating point numbers as a sum of two numbers of the same format. Specifically, given $a$ and $b$ in a specific floating point format, it finds $p$ and $pp$ in the same format so that $p+pp = ab$ exactly. The algorithm from \cite{FPHB} is as follows

\vspace{.2in}
\begin{algorithm} 2MultFMA(a,b)
\hrule
\begin{algorithmic}
   \State{$p = a*b$}
   \State{$pp = fma(a,b,-p)$}
\end{algorithmic}
\hrule
\label{2MultFMA}
\end{algorithm}
\vspace{.2in}

Using this device and then carefully summing will allow us to accurately compute the terms we need. This approach leads us to the following algorithm for a compensated Givens rotation:

\vspace{.2in}
\begin{algorithm} Compensated DLARTG

\hrule
\begin{algorithmic}
   \State{$r = hypot(f,g)$}
   \State{$\bar{c} = f/r$}
   \State{$\bar{s} = g/r$}
   \State{$c_1,c_2 = 2MultFMA(\bar{c},\bar{c})$}
   \State{$s_1,s_2 = 2MultFMA(\bar{s},\bar{s})$}
   \If {$|c| \geq |s|$}
      \State $\epsilon_{Norm} = (1-c_1-s_1-c_2-s_2)/2$
   \Else
      \State $\epsilon_{Norm} = (1-s_1-c_1-s_2-c_2)/2$
   \EndIf
   \State{$p,pp = 2MultFMA(\bar{c},g)$}
   \State{$\epsilon_{Orth} = (fma(-s,f,p)+pp)/r$}
   \State{$\delta_c = \bar{c}*\epsilon_{Norm} - \bar{s}*\epsilon_{Orth}$}
   \State{$\delta_s = \bar{s}*\epsilon_{Norm} + \bar{c}*\epsilon_{Orth}$}
   \State{$\bar{c} = \bar{c} + \delta_c$}
   \State{$\bar{s} = \bar{s} + \delta_s$}
\end{algorithmic}
\hrule
\label{CompensatedDLARTG}
\end{algorithm}
\vspace{.2in}

If no {\tt fma()} is available then we can accomplish the same thing using four calls of another double length product formulation, such the one given in \cite{Dekker}. For instance, in a double precision IEEE 754 format we can use the following algorithm for a double length product:

\vspace{.2in}
\begin{algorithm} 2MultDekker(a,b)
\hrule
\begin{algorithmic}
   \State{split = 134217729.0} \Comment{ A format dependent constant.}
   \State{p = split*a}
   \State{q = a-p}
   \State{HeadofA = q+p}
   \State{TailofA = a-HeadofA}
   \State{p = split*b}
   \State{q = b-p}
   \State{HeadofB = q+p}
   \State{TailofB = b-HeadofB}
   \State{p = a*b}
   \State{pp = (HeadofA*HeadofB-p)+HeadofA*TailofB + HeadofB*TailofA + TailofA*TailofB}
\end{algorithmic}
\hrule
\label{2MultDekker}
\end{algorithm}
\vspace{.2in}

This algorithm can be simplified for computing a square as follows:

\vspace{.2in}
\begin{algorithm} 2SquareDekker(a)
\hrule
\begin{algorithmic}
   \State{split = 134217729.0} \Comment{ A format dependent constant.}
   \State{p = split*a}
   \State{q = a-p}
   \State{HeadofA = q+p}
   \State{TailofA = a-HeadofA}
   \State{p = a*a}
   \State{pp = (HeadofA*HeadofA-p)+2*HeadofA*TailofA + TailofA*TailofA}
\end{algorithmic}
\hrule
\label{2SquareDekker}
\end{algorithm}
\vspace{.2in}

And these can be used to compute the compensated Givens rotation as follows:

\vspace{.2in}
\begin{algorithm} Compensated DLARTG (without the {\tt fma()})

\hrule
\begin{algorithmic}
   \State{$r = hypot(f,g)$}
   \State{$\bar{c} = f/r$}
   \State{$\bar{s} = g/r$}
   \State{$c_1,c_2 = 2SquareDekker(\bar{c})$}
   \State{$s_1,s_2 = 2SquareDekker(\bar{s})$}
   \If {$|c| \geq |s|$}
      \State $\epsilon_{Norm} = (1-c_1-s_1-c_2-s_2)/2$
   \Else
      \State $\epsilon_{Norm} = (1-s_1-c_1-s_2-c_2)/2$
   \EndIf
   \State{$p,pp = 2MultDekker(\bar{c},g)$}
   \State{$q,qq = 2MultDekker(-\bar{s},f)$}
   \State{$\epsilon_{Orth} = (p+q+pp+qq)/r$}
   \State{$\delta_c = \bar{c}*\epsilon_{Norm} - \bar{s}*\epsilon_{Orth}$}
   \State{$\delta_s = \bar{s}*\epsilon_{Norm} + \bar{c}*\epsilon_{Orth}$}
   \State{$\bar{c} = \bar{c} + \delta_c$}
   \State{$\bar{s} = \bar{s} + \delta_s$}
\end{algorithmic}
\hrule
\label{CompensatedDLARTGNoFMA}
\end{algorithm}
\vspace{.2in}

We note that when there is no {\tt fma()} the computational costs rise significantly.

Running the same test as we did in the previous section yields:

\begin{table}[h]
  \begin{center}
      \caption{Error rate for computing the Givens rotation with the compensated algorithm using different variations of the {\tt hypot()} calculation}
      \label{tab:tableGivens}
      \begin{tabular}{l|l|l|l|l|l|l}
         & \multicolumn{2}{c}{Correct} & \multicolumn{2}{c}{Naive} & \multicolumn{2}{c}{Weak}\\
         & Cosine & Sine & Cosine & Sine & Cosine & Sine \\
         \hline
         Zero ulp errors & 100 & 100 & 100 & 100 & 100 & 100 \\
         One ulp errors  & 0 & 0 & 0 & 0 & 0 & 0 \\
         Two ulp errors  & 0 & 0 & 0 & 0 & 0 & 0
      \end{tabular}
   \end{center}
\end{table}

 It is clear that all three approaches lead to perfect observed accuracy for the compensated algorithm and this highly recommends using the compensated approach with the weak hypoteneuse calculation.

\appendix{Appendix A - DLARTG from LAPACK 3.12.0}

This is the core of the DLARTG code from LAPACK 3.12.0 written in FORTRAN. We have excluded the preliminaries to just show the functional core.
\vspace{.2in}

\hrule
\begin{verbatim}
 01  f1 = abs( f )
 02  g1 = abs( g )
 03  if( g == zero ) then
 04     c = one
 05     s = zero
 06     r = f
 07  else if( f == zero ) then
 08     c = zero
 09     s = sign( one, g )
 10     r = g1
 11  else if( f1 > rtmin .and. f1 < rtmax .and. &
 12           g1 > rtmin .and. g1 < rtmax ) then
 13     d = sqrt( f*f + g*g )
 14     c = f1 / d
 15     r = sign( d, f )
 16     s = g / r
 17  else
 18     u = min( safmax, max( safmin, f1, g1 ) )
 19     fs = f / u
 20     gs = g / u
 21     d = sqrt( fs*fs + gs*gs )
 22     c = abs( fs ) / d
 23     r = sign( d, f )
 24     s = gs / r
 25     r = r*u
 26  end if
 27  return
\end{verbatim}
\hrule

% \bibliography{givensbib}{}

\begin{thebibliography}{1}

   \bibitem{Anderson}
   {\sc Anderson, E.}
   \newblock Discontinuous plane rotations and the symmetric eigenvalue problem,
     2001.

   \bibitem{BindelDemmelKahan}
   {\sc Bindel, D., Demmel, J., Kahan, W., and Marques, O.}
   \newblock On computing {G}ivens rotations reliably and efficiently.
   \newblock {\em ACM Trans. Math. Softw. 28}, 2 (June 2002), 206–238.

   \bibitem{borges}
   {\sc Borges, C.~F.}
   \newblock Algorithm 1014: An improved algorithm for hypot(x,y).
   \newblock {\em ACM Trans. Math. Softw. 47}, 1 (Dec. 2020).

   \bibitem{borges2021fast}
   {\sc Borges, C.~F.}
   \newblock Fast compensated algorithms for the reciprocal square root, the
     reciprocal hypotenuse, and {G}ivens rotations.
   \newblock {\em ArXiv e-prints\/} (June 2021).

   \bibitem{Dekker}
   {\sc Dekker, T.~J.}
   \newblock A floating-point technique for extending the available precision.
   \newblock {\em Numerische Mathematik 18}, 3 (1971), 224--242.

   \bibitem{Demmel1997}
   {\sc Demmel, J.}
   \newblock {\em Applied Numerical Linear Algebra}.
   \newblock Other Titles in Applied Mathematics. Society for Industrial and
     Applied Mathematics (SIAM, 3600 Market Street, Floor 6, Philadelphia, PA
     19104), 1997.

   \bibitem{GVL}
   {\sc Golub, G., and Van~Loan, C.}
   \newblock {\em Matrix Computations}.
   \newblock Johns Hopkins Studies in the Mathematical Sciences. Johns Hopkins
     University Press, 2013.

   \bibitem{FPHB}
   {\sc Muller, J., Brunie, N., de~Dinechin, F., Jeannerod, C., Joldes, M.,
     Lef{\`e}vre, V., Melquiond, G., Revol, N., and Torres, S.}
   \newblock {\em Handbook of Floating-Point Arithmetic}.
   \newblock Springer International Publishing, 2018.

   \bibitem{Wilkinson65}
   {\sc Wilkinson, J.}
   \newblock {\em The Algebraic Eigenvalue Problem}.
   \newblock Monographs on numerical analysis. Clarendon Press, 1988.

   \end{thebibliography}
% \bibliographystyle{acm}

\end{document}